\documentclass [12pt,a4paper] {article}
\usepackage{amsfonts}
\usepackage{amsmath}
\usepackage{mathrsfs}
\usepackage{graphicx}
\usepackage{amsmath,amsthm,amssymb,amscd}

\newtheorem{Thm}{Theorem}[section]
\newtheorem{Lem}[Thm]{Lemma}

\theoremstyle{remark}

\theoremstyle{claim}

\DeclareMathOperator{\Diff}{Diff}

\begin{document}

\begin{center}
{\Large \bf  The Structure on Invariant Measures of $C^1$ generic diffeomorphisms}\\
\smallskip
\end{center}
\bigskip
\begin{center}
Wenxiang Sun $^*$
\end{center}
\begin{center}
LMAM, School of Mathematical Sciences, Peking University, Beijing 100871, China\\
\end{center}
\begin{center}
E-mail: sunwx@math.pku.edu.cn
\end{center}
\smallskip
\begin{center}
Xueting Tian $^\dagger$
\end{center}
\begin{center} School of Mathematical Sciences, Peking University, Beijing 100871, China\\
\end{center}
\begin{center}
E-mail: txt@pku.edu.cn
\end{center}
\bigskip

\footnotetext {$^*$ Sun is supported by National Natural Science
Foundation ( \# 10671006, \# 10831003) and National Basic Research
Program of China(973 Program)(\# 2006CB805903) } \footnotetext
{$^\dagger$ Tian is supported by National Natural Science
Foundation( \# 10671006).}
 \footnotetext{AMS Review: 37A25, 37B20, 37C50, 37D20, 37D30;  }

\def\abstractname{\textbf{Abstract}}

\begin{abstract}\addcontentsline{toc}{section}{\bf{English Abstract}}
Let $\Lambda$ be an isolated non-trival transitive set of a $C^1$
generic diffeomorphism $f\in\Diff(M)$. We show that the space of
invariant measures supported on $\Lambda$ coincides with the space
 of accumulation measures of time averages on
one orbit. Moreover, the set of points having this property is
residual in $\Lambda$ (which implies the set of irregular$^+$ points
is also residual in $\Lambda$). As an application, we show that the
non-uniform hyperbolicity of irregular$^+$ points in $\Lambda$ with
totally 0 measure (resp., the non-uniform hyperbolicity of a generic
subset in $\Lambda$) determines the uniform hyperbolicity of
$\Lambda$.
\end{abstract}

\section{Introduction} \setlength{\parindent}{2em}

Let $M$ be a closed $C^{\infty}$ manifold and let $\Diff(M)$ be the
space of diffeomorphisms  of $M$ endowed with the $C^1-$topology.
Denote by $d$ the distance on $M$ induced from a Riemannian metric
on the tangent bundle $TM$. Let $f\in \Diff(M)$.

For a given compact invariant set $\Lambda$, let $P(f|_\Lambda)$ be
the set of periodic points of $f$ in $\Lambda$. Given two periodic
points $p,\,q\in P(f|_\Lambda)$, we say $p,q$ have the barycenter
property, if for any $\varepsilon>0$ there exists an integer
$N=N(\varepsilon,p,q)>0$ such that for any two integers $n_1, n_2$,
there exists a point $x\in P(f|_\Lambda)$ such that $d(f^{i}(z),
f^i(p))<\varepsilon,\,\,-n_1\leq i\leq 0, $ and $d(f^{i+N}(z),
f^i(q))<\varepsilon,\,\,0\leq i\leq n_2$.  $\Lambda$ satisfies the
barycenter property  if the barycenter property holds for any two
periodic points $p,\,q\in P(f|_\Lambda)$.

Given a compact f-invariant set $\Lambda$, $\Lambda$ is transitive
if there is some $x\in\Lambda$ whose forward orbit is dense in
$\Lambda$. A transitive set $\Lambda$ is trivial if it consists of a
periodic orbit. $\Lambda$ is isolated if there is some neighborhood
U of $\Lambda$ in M such that $\Lambda= \cap_{k\in \mathbb{Z}}
f^k(U)$. Denote by ${\cal M}_f(\Lambda),$ ${\cal M}_{erg}(\Lambda)$
and ${\cal M}_p(\Lambda) $ the sets of all $f$ invariant measures,
ergodic measures and periodic measures supported on $\Lambda$
respectively. Clearly ${\cal M}_p(\Lambda) \subseteq {\cal
M}_{erg}(\Lambda) \subseteq {\cal M}_f(\Lambda) $. Given a measure
$\mu$ and $x\in M,$ denote by $V_f(\mu)$ and $V_f(x)$ respectively
the set of all accumulation points of time average
$\mu^N=\frac1N\sum_{j=0}^{N-1}f^j\mu$ and
$\mu^N=\frac1N\sum_{j=0}^{N-1}\delta_{f^jx}$. Note that $V_f(\mu)$
and $V_f(x)$ is a nonempty closed and connected subset of invariant
measures.

Now we state the first theorem as follows.

\begin{Thm}\label{Structure-Thm1}Let $\Lambda$ be an isolated non-trival transitive set of
a $C^1$ generic diffeomorphism $f\in\Diff(M)$. Then the space of
invariant measures supported on $\Lambda$ coincides with the space
of the approximation measures along one orbit, i.e., there is $x\in
\Lambda$ such  that ${\cal M}_f(\Lambda)=V_f(x).$ Moreover, the set
of such points is residual in $\Lambda.$
\end{Thm}

Let us recall the definition of irregular$^+$ point. A point $x\in
M$ is called irregular for positive iterations (or shortly
irregular$^+$) if there is a continuous function $\phi: M
\rightarrow\mathbb{R}$ such that the sequence
$\frac1n\sum_{i=0}^{n-1}\phi(f^i(x))$ is not convergent. Cleary
every point $x$ with ${\cal M}_f(\Lambda)=V_f(x)$ is irregular$^+$.
By Theorem \ref{Structure-Thm1}, these points are "many" for generic
diffeomorphisms. But by Birkhoff ergodic theorem, the set of
irregular$^{+}$ points is a totally probability 0 measure set, i.e.,
for any invariant measure, its measure is zero. So irregular$^{+}$
points is "few" in the probabilistic perspective. It is a very
interesting phenomena.

We recall  the notions of uniform hyperbolicity and non-uniform
hyperbolicity. Let $f :M\rightarrow M$ be a diffeomorphism on a
compact manifold M. A compact invariant set $\Delta$ of $f$ is
called hyperbolic if there is a continuous invariant splitting
$T_\Delta M = E\oplus
 F$  and two constants
 $C > 0, 0 < \lambda < 1,$ such
that $$\|Df^ n|_{E(x)}\|\leq C\lambda^n,\,\,\,and \,\,\,\|Df^
{-n}|_{F(x)}\|\leq C\lambda^n, \forall n \in \mathbb{N},\,\,\forall
x\in \Delta.$$ We say that a point $x\in M$ is a $NUH$ point or,
simply, $NUH$, if\\
 (1) there is a $Df$-invariant splitting $T_{Orb(x)}M = E_{Orb(x)}\oplus
 F_{Orb(x)}$,\\
 (2) there exist two constants $\eta> 0$, $L \in \mathbb{N}$ and a  Riemannian
 metric $\|\cdot\|$ such that
$$\limsup_{n\rightarrow\infty}\frac1n\sum_{j=0}^{n-1}\log\|Df^L(f^j
(x))|_{E (f^j (x))}\|\leq-\eta$$ and
$$\limsup_{n\rightarrow\infty}\frac1n\sum_{j=0}^{n-1}\log\|[Df^L(f^j
(x))|_{F (f^j (x))}]^{-1}\|\leq-\eta.$$

The following Theorem shows that the non-uniform hyperbolicity of a
totally 0 measure set  determines the uniform hyperbolicity of the
whole space. More precisely, for an isolated non-trival transitive
set $\Lambda$ of a generic diffeomorphism, the non-uniform
hyperbolicity of irregular$^+$ points in $\Lambda$ determines the
uniform hyperbolicity of $\Lambda$.

\begin{Thm}\label{Structure-Thm2}Let $\Lambda$ be an isolated non-trival transitive set of
a $C^1$ generic diffeomorphism $f\in\Diff(M)$ and let $T_{\Lambda}
M=E\oplus F$ be a continuous $Df$-invariant splitting over
$\Lambda$. Let $U\subset \Lambda$ be a nonempty open set. If
(generic) irregular$^+$ points in $U$ satisfy $NUH$ condition with
respect to $T_{\Lambda} M=E\oplus F$, then $\Lambda$ is a hyperbolic
basic set for $f.$
\end{Thm}

\section {\bf Proof of our Theorem \ref{Structure-Thm1}}

In this section we suppose the assumptions of following lemmas are
all the same as Theorem \ref{Structure-Thm1}.
\begin{Lem}\label{Structure-Lem5}
(i) The set of periodic measures supported on $\Lambda$  is a dense
subset of the set ${\cal M}_f(\Lambda)$ of invariant measures
supported on $\Lambda$: $\overline{\mathcal{M}_{p}(\Lambda)}={\cal
M}_f(\Lambda).$

 (ii)
$\Lambda$ satisfies the barycenter property.
\end{Lem}

{\bf Proof}   The proof of (i) and (ii) are given in \cite{ABC}, see
Theorem 3.5 and Proposition 4.8, respectively. \hfill $\Box$

An argument by Bonatti and D\'{\i}az \cite{BD}, based on Hayashi
Connecting Lemma\cite{H}, shows that isolated transitive sets
$\Lambda$ of $C^1$ generic diffeomorphisms are relative homoclinic
classes:

\begin{Lem}\label{Structure-Lem3}(\cite{BD}) There
is some periodic point $p$ such that $\Lambda=H(p),$ where $H(p)$
denotes the homoclinic class of $p.$
\end{Lem}

Here we divide into the following two lemmas to prove Theorem
\ref{Structure-Thm1}.
\begin{Lem}\label{Structure-Lem2} There
is $x\in \Lambda$ such  that ${\cal M}_f(\Lambda)=V_f(x).$ Moreover,
the set of such points is dense in $\Lambda.$
\end{Lem}

 {\bf Proof}
Since ${\cal M}_f(\Lambda)$ is closed and connected, there exists a
sequence of closed balls $B_n$ in ${\cal M}_f(\Lambda)$ with radius
$\varepsilon_n$ (in some metric $\widetilde d$ compatible with the
weak$^*$ topology) such that
the following holds: \\
(a) $B_n\cap B_{n+1}\neq \emptyset,$\\
(b) ${\cap_{N=1}^{\infty}\cup_{n\geq N}} B_n={\cal M}_f(\Lambda),$ \\
(c) $\lim_{n\rightarrow+\infty}\varepsilon_n=0$.\\
By Lemma \ref{Structure-Lem5} (i), ${\cal M}_f(\Lambda)=
\overline{{\cal M}_p(\Lambda)}$. We may also assume that the center
of $B_n$ is a periodic measure $Y_n$. The support of $Y_n$ is the
orbit of some periodic point $x_n \in \Lambda$ whose period is
$p_n$.

Let $x_0 \in \Lambda$ be given and $U_0$ the open ball of radius
$\delta$ around $x_0$. By Lemma \ref{Structure-Lem3}, the set of
periodic points in $\Lambda$ is a dense subset. Hence, without loss
of generality, we can assume $x_0$ is a periodic point. Since ${\cal
M}_f(\Lambda)\supseteq V_f(x)$ is trivial, we only need to show that
there exists an x $\in U_0$ such that ${\cal M}_f(\Lambda)\subseteq
V_f(x)$.

{\bf Step 1} choose some x $\in U_0$ that we need.\\
 Let $z_0=x_0.$ By Lemma \ref{Structure-Lem5} (ii), $\Lambda$ satisfies barycenter property. For
periodic points $z_0=x_0$ and $x_1,$ there exists a positive integer
$M_{1}$ and a periodic point $z_1 \in \Lambda$ such that
$$d(f^jz_0,f^jz_1)<2^{-1}\delta\,\,\,\,\,\text{for}\,\,j=a_0=b_0=0$$
and
$$d(f^jx_1,f^jz_1)<2^{-1}\delta\,\,\,\,\,\text{for}\,\,a_1\leq j\leq
b_1,$$ where
$$a_0=0,b_0=0$$
$$a_1=b_0+M_1,b_1=a_1+2(b_0+M_{1})p_1.$$

Using barycenter property again, for periodic points $z_1$ and
$x_2,$ we have a positive integer $M_{2}$ and  a
  periodic point $z_2 \in \Lambda$ such that
$$d(f^jz_1,f^jz_2)<2^{-2}\delta\,\,\,\,\,for\,\,0=a_0\leq j\leq
b_1$$ and
$$d(f^jx_2,f^jz_2)<2^{-2}\delta\,\,\,\,\,for\,\,a_2\leq j\leq
b_2,$$ where
$$a_2=b_1+M_2,b_2=a_2+2^2(b_1+M_{2})p_2.$$

In general  we have a positive integer $M_n $ and  a
 periodic point $z_n \in \Lambda$ such that
$$d(f^jz_{n-1},f^jz_n)<2^{-n}\delta\,\,\,\,\,for\,\,0=a_0\leq j\leq
b_{n-1}$$ and
$$d(f^jx_n,f^jz_n)<2^{-n}\delta\,\,\,\,\,for\,\,a_n\leq j\leq
b_n,\,\,\,\,\,\,$$ where
$$a_n=b_{n-1}+M_n,b_n=a_n+2^n(b_{n-1}+M_{n})p_n.$$

It is easy to check that for $m>n$,
$$
d(f^jx_n,f^jz_m)<2^{-n+1}\delta\,\,\,\,\,for\,\,a_n\leq j\leq
b_n\,.\,\,\,\,\,
$$
Since $d(z_{n-1},z_{n})<2^{-n}\delta,$ the sequence $z_n$ converges
to some point $x\in U_0$, and one has
\begin {equation}\label{Structure-Eq:5}
d(f^jx_n,f^jx)<2^{-n+1}\delta\,\,\,\,\,for\,\,a_n\leq j\leq
b_n\,.\,\,\,\,\,
\end {equation}
Remark that if A is a finite subset of $\mathbb{N},$ then
\begin {equation}\label{Structure-Eq:6}
|\frac1{card A}\sum_{j\in A}\xi(f^jy)-\frac1{max
A+1}\sum_{j=0}^{maxA}\xi(f^jy)|\leq 2(cardA)^{-1}(max A+1-card
A)||\xi||
\end {equation}
 for any $y \in M $ and $\xi \in C(M)$.

{\bf Step 2} ${\cal M}_f(\Lambda)\subseteq V_f(x)$.

Let $\nu \in {\cal M}_f(\Lambda)$ be given. By (b) and (c) there
exists an increasing sequence $n_k\uparrow\rightarrow\infty$ such
that $Y_{n_k}\rightarrow\nu$. Let $\xi\in C(M)$ be given with
$||\xi||\leq1$, and denote by $w_\xi(\varepsilon)$ the oscillation
$$max\{|\xi(y)-\xi(z)|\,\,\big{|}\,\,d(y,z)\leq\varepsilon\}.$$
Let $\nu_n$ denote the measure $\delta(x)^{b_n}$. Thus $$\int \xi d
\nu_n=\frac1{b_n}\sum_{j=0}^{b_n-1}\xi(f^jx).$$
 Also $$\int \xi d
 Y_n=\frac1{b_n-a_n}\sum_{j=a_n}^{b_n-1}\xi(f^jx_n).$$
So by (\ref{Structure-Eq:5}),$$|\int \xi d
 Y_n-\frac1{b_n-a_n}\sum_{j=a_n}^{b_n-1}\xi(f^jx)|\leq w_{\xi}(2^{-n+1}\delta).$$
 Since $||\xi||\leq1$, (\ref{Structure-Eq:6}) implies, with $A=[a_n,b_n]\cap
 Z$,$$|\frac1{b_n-a_n}\sum_{j=a_n}^{b_n-1}\xi(f^jx)-\frac1{b_n}\sum_{j=0}^{b_n-1}\xi(f^jx)|\leq\frac{2a_n}{b_n-a_n}.$$
 Since $w_{\xi}(2^{-n+1}\delta)\rightarrow0$ and $\frac{2a_n}{b_n-a_n}\rightarrow
 0$ as $n\rightarrow+\infty,$ this shows that $$|\int \xi d
\nu_n-\int \xi d
 Y_n|\rightarrow0.$$ Hence $\nu_{n_k}\rightarrow\nu$ and thus $\nu \in V_f(x).$
 \hfill $\Box$

\begin{Lem}\label{Structure-Lem6} The set $\{x\in\Lambda\,|\,{\cal M}_f(\Lambda) =V_f(x)\}$ is residual in $\Lambda.$
\end{Lem}
{\bf Proof}  Denote by ${\cal M}(\Lambda)$ the set of all Borel
probability measures defined on $\Lambda$. Take open balls
$V_i,\,\,U_i$ ($i\in\mathbb{N}$) in ${\cal M}(\Lambda)$ such that \\
(a) $V_i\subseteq\overline{V_i}\subseteq U_i;$ \\
(b) $diam(U_i)\rightarrow0;$\\
(c) $V_i\cap{\cal M}_f(\Lambda)\neq\emptyset;$ \\
(d) each point of ${\cal M}_f(\Lambda)$ lies in infinitely many
$V_i.$\\
Put $$P(U_i)=\{x\in\Lambda\,|\, V_f(x)\cap U_i\neq\emptyset\}.$$ It
is easy to see that the set of points with ${\cal M}_f(\Lambda)
=V_f(x)$ is just $\cap_{i\in\mathbb{N}} P(U_i).$ Since $V_i\cap{\cal
M}_f(\Lambda)\neq\emptyset,$  one has
$$\cap_{i\in\mathbb{N}} P(U_i)=\{x\in\Lambda\,|\,{\cal M}_f(\Lambda) =V_f(x)\}\subseteq\{x\in\Lambda\,|\,\forall N_0\in\mathbb{N},\,
\exists N>N_0 \,\,\,{with}\,\,\, \delta(x)^N\in V_i\}$$
$$=\cap_{N_0=1}^{\infty}\,\cup_{N>N_0}\,\{x\in
\Lambda\,\,|\,\,\delta(x)^N\,\in \,V_i\}.$$ By the definition of
$P(U_i),$ for any $i$, $$
P(U_i)\supseteq\cap_{N_0=1}^{\infty}\,\cup_{N>N_0}\,\{x\in
\Lambda\,\,|\,\,\delta(x)^N\,\in \,V_i\}.$$ So
$$\{x\in\Lambda\,|\,{\cal M}_f(\Lambda) =V_f(x)\}=\cap_{i\in\mathbb{N}}
P(U_i)=\cap_{i\in\mathbb{N}}\cap_{N_0=1}^{\infty}\,\cup_{N>N_0}\,\{x\in
\Lambda\,\,|\,\,\delta(x)^N\,\in \,V_i\}.$$ Since
$$x\mapsto\delta(x)^N$$ is continuous (for fixed $N$), the sets of
$\cup_{N>N_0}\,\{x\in \Lambda\,\,|\,\,\delta(x)^N\,\in \,V_i\}$ are
open. By Lemma \ref{Structure-Lem2}, the sets of
$$\cup_{N>N_0}\,\{x\in \Lambda\,\,|\,\,\delta(x)^N\,\in
\,V_i\}\supseteq\{x\in\Lambda\,|\,{\cal M}_f(\Lambda) =V_f(x)\}$$
are dense in $\Lambda.$ Hence $\{x\in\Lambda\,|\,{\cal M}_f(\Lambda)
=V_f(x)\}$ is residual in $\Lambda.$\hfill $\Box$

\bigskip
\section {\bf Proof of our Theorem \ref{Structure-Thm2}}

Before proving Theorem \ref{Structure-Thm2} we need the following
lemma\cite{Cao} by Y.Cao.

\begin{Lem}\label{Structure-Lem4} Let $f:M\rightarrow M$ be a $C^1$ local diffeomorphism on a compact
manifold and let $\Lambda$ be a compact and $f-$invariant set.
Suppose that there exists a continuous $Df$-invariant splitting
$T_\Lambda M=E\oplus F$. If the Lyapunov exponents restricted on $E$
and $F$ of every $f$ invariant probability measure are all negative
and positive respectively, then $\Lambda$ is uniformly hyperbolic.
\end{Lem}

The following lemma shows that the NUH condition of a point $x$ with
${\cal M}_f(\Lambda)= V_f(x)$ determines the uniform hyperbolicity
of $\Lambda$, which can deduce Theorem \ref{Structure-Thm2}.

\begin{Lem}\label{Structure-Lem1} Let $f:M\rightarrow M$ be a $C^1$ local diffeomorphism on a compact
manifold and let $\Lambda$ be a compact and $f-$invariant set.
Suppose that there exists a continuous $Df$-invariant splitting
$T_\Lambda M=E\oplus F$. If $y$ is a NUH point in $\Lambda$ with
respect to $T_y M=E_y\oplus F_y$ and ${\cal M}_f(\Lambda)= V_f(y)$,
then $\Lambda$ is uniformly hyperbolic basic set.
\end{Lem}

{\bf Proof}  By assumption, take $\eta> 0$, $L \in \mathbb{N}$ and a
Riemannian
 metric $\|\cdot\|$ such that
$$\limsup_{n\rightarrow\infty}\frac1n\sum_{j=0}^{n-1}\log\|Df^L(f^j
(x))|_{E (f^j (y))}\|\leq-\eta$$ and
$$\limsup_{n\rightarrow\infty}\frac1n\sum_{j=0}^{n-1}\log\|[Df^L(f^j
(x))|_{F (f^j (y))}]^{-1}\|\leq-\eta.$$ Let
$$\varphi_E(x)=log\|Df^L|_{E(x)}\|,x\in \Lambda.$$ By the continuity
of $T_\Lambda M=E\oplus F$,  $\varphi_E(x)$ is continuous on
$\Lambda.$ Since ${\cal M}_f(\Lambda)= V_f(y)$, for any given
$\mu\in{\cal M}_f(\Lambda)$, there is $n_k\uparrow+\infty$ such that
$$\frac1{n_k}\sum_{i=0}^{n_k}\delta_{f^i(y)}\rightarrow \mu$$ in the
weak* topology. So
$$\int\varphi_E(x)d\mu
=\lim_{n\rightarrow+\infty}\frac1{n_k}\sum_{i=0}^{n_k}\varphi_E{f^i(y)}\leq-\eta.
$$

 By Birkhorff Ergodic Theorem,
$$\lim_{n\rightarrow+\infty}\frac1n\sum_{i=0}^{n-1}\varphi_E(f^ix)$$ exists on a totally measure set.
We claim that there exists a totally measured set $\Delta$ such that
$$\lim_{n\rightarrow+\infty}\frac1n\sum_{i=0}^{n-1}\varphi_E(f^ix)\leq-\eta$$
for any $x\in\Delta.$ Otherwise, there exists an invariant measure
$\mu$ and a $\mu-$positive measure set $\Gamma$ such that
$$\lim_{n\rightarrow+\infty}\frac1n\sum_{i=0}^{n-1}\varphi_E(f^ix)>-\eta,$$ for any
$x\in\Gamma.$ Since the  limit function above is $f-$invariant, we
can assume that $\Gamma$ is $f-$invariant. So we can define an
invariant measure $\nu$  as follows: for every Borel set $B$,
$$\nu(B):=\frac{\mu(B\cap\Gamma)}{\mu(\Gamma)}.$$ Then $\nu$ is
an invariant measure and $\nu(\Gamma)=1$. So by Birkhorff Ergodic
Theorem, we have
$$\int\varphi_E(x)d\nu=\int\lim_{n\rightarrow+\infty}\frac1n\sum_{i=0}^{n-1}\varphi_E(f^ix)d\nu
>-\eta,$$ which contradicts $\int\varphi_E(x)d\nu\leq-\eta.$

By Oseledec theorem, the limit
$$\lambda_E(x):=\lim_{n\rightarrow+\infty}\frac1n\log\|Df^n|_{E(x)}\|$$
exists on a totally measured set $\Delta'$ (i.e., $\lambda_E(x)$ is
the maximal Lyapunov exponent of $x$ on subbundle $E(x)$). Since
$\lambda_E(x)$ is $f-$invariant, then by sub-addition of
$\log\|Df^n|_E(x)\|$ we have
\begin{eqnarray*}&&\lambda_E(x)=\frac1L\sum_{j=0}^{L-1}\lambda_E(f^jx)\\
=&&\lim_{n\rightarrow+\infty}\frac1L\sum_{j=0}^{L-1}\frac1{nL}\log\|Df^{nL}|_{E(f^jx)}\|\\
\leq&&\lim_{n\rightarrow+\infty}\frac1L\sum_{j=0}^{L-1}\frac1{nL}\sum_{i=0}^{n-1}\log\|Df^{L}|_{E(f^{j+iL}x)}\|\\
=&&\frac1L\lim_{n\rightarrow+\infty}\frac1{nL}\sum_{i=0}^{nL-1}\varphi_E(f^ix)\\
\leq&&-\frac{\eta}L<0
\end{eqnarray*}
for all $x\in\Delta\cap\Delta'$. By Lemma \ref{Structure-Lem4} one
gets that $E$ is a contracting subbundle. Similarly we also have
that $F$ is an expanding subbundle.
\hfill $\Box$\\

\end{document}